\documentclass[a4paper,10pt]{article}

\pdfoutput=1

\usepackage[paper=letterpaper,margin=1in]{geometry}

\usepackage{amsmath,amssymb,amsfonts}

\begin{document}
\newcommand{\todo}[1]{{\bf ?????!!!! #1 ?????!!!!}\marginpar{$\Longleftarrow$}}
\newcommand{\fref}[1]{Figure~\ref{#1}}
\newcommand{\tref}[1]{Table~\ref{#1}}
\newcommand{\sref}[1]{\S~\ref{#1}}
\newcommand{\nn}{\nonumber}
\newcommand{\tr}{\mathop{\rm Tr}}
\newcommand{\ch}{\rm Ch}
\newcommand{\comment}[1]{}

\newcommand{\cM}{{\cal M}}
\newcommand{\cW}{{\cal W}}
\newcommand{\cN}{{\cal N}}
\newcommand{\cH}{{\cal H}}
\newcommand{\cK}{{\cal K}}
\newcommand{\cZ}{{\cal Z}}
\newcommand{\cO}{{\cal O}}
\newcommand{\cB}{{\cal B}}
\newcommand{\cC}{{\cal C}}
\newcommand{\cD}{{\cal D}}
\newcommand{\cE}{{\cal E}}
\newcommand{\cF}{{\cal F}}
\newcommand{\cR}{{\cal R}}
\newcommand{\IA}{\mathbb{A}}
\newcommand{\IP}{\mathbb{P}}
\newcommand{\IQ}{\mathbb{Q}}
\newcommand{\IH}{\mathbb{H}}
\newcommand{\IR}{\mathbb{R}}
\newcommand{\IC}{\mathbb{C}}
\newcommand{\IF}{\mathbb{F}}
\newcommand{\IM}{\mathbb{M}}
\newcommand{\II}{\mathbb{I}}
\newcommand{\IZ}{\mathbb{Z}}
\newcommand{\re}{{\rm Re}}
\newcommand{\im}{{\rm Im}}
\newcommand{\sym}{{\rm Sym}}

\newcommand{\tmat}[1]{{\tiny \left(\begin{matrix} #1 \end{matrix}\right)}}
\newcommand{\mat}[1]{\left(\begin{matrix} #1 \end{matrix}\right)}
\newcommand{\diff}[2]{\frac{\partial #1}{\partial #2}}
\newcommand{\gen}[1]{\langle #1 \rangle}
\newcommand{\ket}[1]{| #1 \rangle}
\newcommand{\jacobi}[2]{\left(\frac{#1}{#2}\right)}

\newcommand{\drawsquare}[2]{\hbox{%
\rule{#2pt}{#1pt}\hskip-#2pt
\rule{#1pt}{#2pt}\hskip-#1pt
\rule[#1pt]{#1pt}{#2pt}}\rule[#1pt]{#2pt}{#2pt}\hskip-#2pt
\rule{#2pt}{#1pt}}
\newcommand{\fund}{\raisebox{-.5pt}{\drawsquare{6.5}{0.4}}}
\newcommand{\antifund}{\overline{\fund}}

\renewcommand{\thefootnote}{\fnsymbol{footnote}}

~\\

\begin{center}
{\Large \bf Moonshine and the Meaning of Life}
\end{center}

\centerline{
{\large Yang-Hui He}$^1$ \&
{\large John McKay}$^2$
}
\vspace*{1.0ex}

{\it
{\scriptsize
\begin{tabular}{ll}
\begin{tabular}{c}
${}^{1}$ 
Dept.~of Maths, City U., London, EC1V 0HB, UK; \\
School of Physics, NanKai U., Tianjin, 300071, P.R.~China; \\
Merton College, University of Oxford, OX14JD, UK\\
\qquad hey@maths.ox.ac.uk\\
\end{tabular}
&
\begin{tabular}{c}
${}^{2}$ 
Department of Mathematics and Statistics,\\
Concordia University, 1455 de Maisonneuve West,\\
Montreal, Quebec, H3G 1M8, Canada\\
\qquad mckay@encs.concordia.ca
\end{tabular}
\end{tabular}
}}


\comment{
\begin{abstract}
With a jocund air, we present an observation on the first 24 coefficients of the modular invariant and of the modular discriminant.
The observation is purely for the sake of entertainment and could be of some diversion to a mathematical audience.
\end{abstract}
}

\setcounter{equation}{0}

\thispagestyle{empty}

\vspace{.5cm}

The elliptic modular function, $j$, invariant under $PSL(2,\IZ)$, has Fourier expansion 
\begin{equation}
j(q) = \frac{E_4(q)^3}{\Delta(q)} = \sum\limits_{m=-1}^\infty c_m q^m = 
\frac{1}{q} + 744 + 196884q + 21493760q^2 + \ldots \ , 
\end{equation}
as $z \to i\infty$, where $q = e^{2 \pi i z}$ is the nome for $z$.
$E_4(z) 
= 1 + 240 \sum\limits_{n=1}^\infty \sigma_3(n) q^n$ is the theta series for the $E_8$ lattice, $\sigma_3(n) = \sum\limits_{d|n} d^3$ and
\begin{equation}\label{delta}
\Delta(q) = q \prod \limits_{n=1}^\infty (1 - q^n)^{24} = \sum\limits_{m=1}^\infty \tau_m q^m 
= q -24q^2+252q^3-1472q^4+4830q^5 + \ldots 
\end{equation}
is the modular discriminant \cite{serre}.
There are two new congruences
\begin{quote}
{\bf OBSERVATIONS:} \quad
$\bullet$ [JM] 
$\left( \sum\limits_{m=1}^{24} c_m^2 \right) \bmod 70 \equiv 42$ \ ; \quad
$\bullet$ [YHH] 
$\left( \sum\limits_{m=1}^{24} \tau_m^2 \right) \bmod 70 \equiv 42$ \ .
\end{quote}

The vector $\omega = (0,1,2,\ldots,24 : 70)$ lives in the Lorentzian lattice $II_{25,1}$ in 26 dimensions as an isotropic Weyl vector \cite{con}, allowing us to construct 
\comment{
\footnote{
Indeed, the Diophantine identity
$1^2 + 2^2 + \ldots + 24^2 = 70^2$
is the only one of its kind, involving an integer quadratic form in as high as 24 variables.
In fact, 24 is the only integer greater than 1 such that consecutive squares summing up to which gives itself a perfect square \cite{lucas,watson}. 
}}
the Leech lattice
as $\omega^\perp / \omega$. 
Watson's \cite{lucas,watson} unique non-trivial solution to $\sum\limits_{i=1}^{n} i^2 = m^2$ is $(n,m)=(24,70)$.

Indeed, the second author's observation 35 years ago that
\begin{equation}
196884 = 196883 + 1 
\end{equation}
sparked the field of ``Monstrous Moonshine'' \cite{borcherds,CN}, underlying so much mathematics and physics, relating, inter alia, modular functions, finite groups, lattices, conformal field theory, string theory and gravity (see \cite{Gannon:2004xi} for a review of some of the vast subjects encompassed) in which the $j$-invariant and the Leech lattice are central.
As we ponder the meaning of life, we should be aware of the prescient remarks of the author \cite{adams}, Douglas Adams:
\begin{quote}
{\it ``The Answer to the Great Question \ldots is \ldots Forty-two,'' said Deep Thought, with infinite majesty and calm.}
\end{quote}

\end{document}